\documentclass[11pt, oneside]{article}
\usepackage{amsfonts}
\usepackage{mathrsfs}
\usepackage{color}
\usepackage[colorlinks]{hyperref}
\usepackage{latexsym}
\usepackage{amssymb}
\usepackage{amsmath}
\usepackage{enumerate}
\usepackage{amsthm}
\usepackage{indentfirst}
\usepackage{mathtools}
\usepackage{tikz}
\usepackage{psfrag}
\usepackage{graphicx}
\usepackage{bm}
\usepackage[rflt]{floatflt}
\usepackage{float}
\usepackage{authblk}
\usepackage[OT2,T1]{fontenc}
\usepackage{ulem}
\usepackage[square, comma, sort&compress, numbers]{natbib}
\usepackage{verbatim}

\numberwithin{equation}{section}
\allowdisplaybreaks

\newtheorem{theorem}{\color{black}\indent Theorem}[section]
\newtheorem{lemma}{\color{black}\indent Lemma}[section]
\newtheorem{proposition}{\color{black}\indent Proposition}[section]
\newtheorem{definition}{\color{black}\indent Definition}[section]
\newtheorem{remark}{\color{black}\indent Remark}[section]

\usepackage{amssymb,amsmath}
\begin{document}
\title{\LARGE\bf
A relativistic quantum Euler–Poisson system derived from the Klein–Gordon–Poisson equation: hyperbolic–elliptic structure}

\author{
Ben Duan$^{1}$,
Bin Guo$^{1}$,
Jun Li$^{2}$,
Rongrong Yan$^{1}$\thanks{Corresponding author. Email addresses:
bduan@jlu.edu.cn (Ben Duan),
bguo@jlu.edu.cn (Bin Guo),
lijun@nju.edu.cn (Jun Li),
yanrr2024@gmail.com (Rongrong Yan)}
}

\date{}
\maketitle
\vspace{-2.3em}

\begin{center}
$^{1}$School of Mathematics, Jilin University, Changchun 130012, PR China\\
$^{2}$School of Mathematics, Nanjing University, Nanjing 210093, PR China
\end{center}
\abstract{In the Klein–Gordon equation, quantum and relativistic parameters are strongly coupled,
which poses significant analytical challenges in the derivation and analysis of related
classical fluid models. In this paper, starting from the Klein–Gordon–Poisson system,
we formally derive a relativistic quantum hydrodynamic (RQHD) system via the Madelung
transformation, in which the relativistic and quantum correction terms in the Euler–Poisson framework are clearly exhibited. In particular, at a formal level, the RQHD system reduces to the relativistic hydrodynamics system in the semiclassical regime and to the quantum hydrodynamics system in the non-relativistic regime.  These limiting
procedures highlight the unified structure of the proposed model and clarify the role played by the coupled relativistic and quantum effects. From an analytical point of view, by reformulating the RQHD system as a coupled
hyperbolic–elliptic system with a nonlocal Poisson interaction, we establish the
local-in-time existence and uniqueness of classical solutions to the associated Cauchy
problem. The initial density is assumed to be a small perturbation of a positive constant
state, while the remaining initial data are taken to be general smooth functions. The
analysis relies on energy estimates and suitable estimates for the nonlocal terms,
and provides a rigorous well-posedness result in the natural energy space.}\\
\\
\maketitle\textbf{Keywords}: Klein--Gordon--Poisson system; relativistic quantum hydrodynamic system; Euler–Poisson system; Parameter Limits; hyperbolic--elliptic system; local well-posedness.
 \thispagestyle{empty}
\section{Introduction}

This paper is concerned with a relativistic quantum hydrodynamic model motivated by the self-consistent Klein--Gordon--Poisson system. As a starting point, we consider the following Klein--Gordon equation coupled with a Poisson equation:
\begin{equation}\label{AA}
\frac{\hbar^2}{2mc^2}\partial_t^2\varphi
-\frac{\hbar^2}{2m}\Delta\varphi
+\frac{mc^2}{2}\varphi
+V(x,t)\varphi=0,
\end{equation}
posed on the three-dimensional space $\mathbb{R}^3$.
Here $m>0$ denotes the particle mass, $c$ is the speed of light, and $\hbar$ is the reduced Planck constant. The unknown $\varphi=\varphi(x,t)$ is a complex-valued scalar field defined on $\mathbb{R}^{3+1}$. The potential $V$ is determined by the Poisson equation
\[
-\Delta V = |\varphi|^2 - b(x) \quad \text{in } \mathbb{R}^3,
\]
where $b(x)$ denotes a prescribed background charge density.

A classical question associated with the Klein--Gordon equation \eqref{AA} concerns its non-relativistic regime, which formally corresponds to the limit where the speed of light $c$ becomes large. In this regime, it is natural to ask whether solutions of \eqref{AA} can be approximated by solutions of the Schr\"odinger equation
\begin{equation}\label{Hhh1}
i\hbar \partial_t \phi + \frac{\hbar^2}{2m}\Delta \phi - V(x,t)\phi = 0.
\end{equation}
Such non-relativistic approximations have been studied in the literature \cite{wz bao,RS3,RS4}. Under suitable assumptions, these works show that solutions of Klein--Gordon-type equations admit Schr\"odinger-type approximations as $c\to\infty$.  In \cite{SML2}, Zhang established the convergence from the Schrödinger–Poisson system to the Vlasov–Poisson system in one spatial dimension. Most existing studies concentrate on the semiclassical limit of the Schrödinger–Poisson system; see, for example, \cite{SML4,L.S. Mai,SML3, Q.C. Ju}. In addition, \cite{G.J. Zhang} developed a semiclassical limit theory for the bipolar quantum hydrodynamic model arising in semiconductor physics.

In contrast, rigorous analytical results on the semiclassical limit of the nonlinear Klein--Gordon equation remain relatively scarce. In \cite{inF1}, the authors investigate the convergence of weak solutions with finite charge energy to relativistic wave maps in the semiclassical limit. On the other hand, \cite{inF3} studies the convergence of the energy associated with the nonlinear Klein--Gordon equation toward the energy of the compressible Euler equations in a non-relativistic--semiclassical limit.

Due to the coupling between the quantum and relativistic parameters in the Klein–Gordon equation, a direct analysis of the semiclassical or non-relativistic–semiclassical limits at the level of the original model is not straightforward. Modulated energy methods have been successfully applied to describe the asymptotic behavior of solutions in the non-relativistic regime. However, their application to the semiclassical limit remains difficult.

Motivated by this gap, we apply the Madelung transformation \cite{MdL} to reformulate the Klein--Gordon--Poisson equation as a hydrodynamic system incorporating both relativistic and quantum effects. This reformulation leads to the following relativistic quantum hydrodynamic (RQHD) system:
\begin{align}\label{J1}
\begin{cases}
\partial_t n + \mathrm{div}(n\nabla S) = \upsilon^2 \partial_t (n \partial_t S), \\[2mm]
\partial_t (n\nabla S)
+ \mathrm{div}\!\left(\dfrac{n\nabla S\otimes n\nabla S}{n}\right)
- \dfrac{\varepsilon^2}{2} n \nabla\!\left(\dfrac{\Delta \sqrt{n}}{\sqrt{n}}\right)
+ n \nabla V \\[2mm]
\qquad\qquad
= \dfrac{\upsilon^2}{2}
\Bigl[ 2 \partial_t (\partial_t S \nabla S\, n)
- \varepsilon^2 \partial_t \bigl(n \nabla (\partial_t \log n)\bigr) \Bigr], \\[2mm]
-\Delta V = n-b(x) .
\end{cases}
\end{align}
Here $n$ and $S$ denote the particle density and the phase, respectively, while $\varepsilon>0$ and $\upsilon>0$ denote the quantum and relativistic parameters. The relativistic quantum hydrodynamic system offers a formulation in which the formal reductions corresponding to the non-relativistic and semiclassical limits can be observed at a formal level in a transparent manner. In this context, the well-posedness analysis of the associated hydrodynamic system may be regarded as a preliminary step toward the study of parameter-dependent limits.

Another objective of this paper is to establish the local-in-time existence and uniqueness of classical solutions to the relativistic quantum hydrodynamic system. No smallness assumptions are imposed on the initial data, except for the density.

Guo \cite{Y. Guo} addressed the Poisson coupling in the Euler–Poisson system by introducing a Riesz potential representation for the electrostatic potential, thereby rewriting the system as an evolution equation with nonlocal terms. This approach reveals the intrinsic hyperbolic–elliptic structure of the system and lays the analytical foundation for its well-posedness theory. Subsequently, Gu and Lei \cite{X. Gu} employed the Newtonian potential representation in their study of the Euler–Poisson equations with physical vacuum boundaries, and established the local well-posedness of the system within this framework. On the other hand, for Euler–Poisson systems with quantum corrections, Li and Marcati \cite{H.L. Li} proved the local-in-time existence of classical solutions to multidimensional quantum hydrodynamic models. Subsequently, Li et al. \cite{C.C. Hao} extended this result to the three-dimensional whole space $\mathbb{R}^3$ and further \cite{H.L. Li2} established the global-in-time existence and uniqueness of strong solutions.

In the present work, we develop a hyperbolic–elliptic analytical framework tailored to the coupled system under consideration: the electrostatic potential is governed by a Poisson equation, while the remaining hydrodynamic variables form a hyperbolic subsystem. This formulation yields a mathematically tractable structure provides the analytical framework adopted in this paper for the local well-posedness analysis. 

\medskip
\noindent\textbf{Formal singular limits and connections with classical models.}
We briefly discuss several singular limits of the relativistic quantum hydrodynamic system at a purely formal level, without attempting any rigorous justification, in order to illustrate its connections with classical fluid models.

\medskip
\noindent\textbf{Formal semiclassical limit.}
To the best of our knowledge, rigorous results on the semiclassical limit of the nonlinear Klein--Gordon equation remain rather limited.
At the level of the relativistic quantum hydrodynamic system \eqref{J1}, as quantum effects vanish $(\varepsilon \to 0)$, the system formally reduces to the relativistic Euler--Poisson system
\begin{equation}\label{EP-rel}
\begin{cases}
\partial_t n + \mathrm{div}(n\nabla S) = \upsilon^2 \partial_t (n \partial_t S),\\[2mm]
\partial_t (n\nabla S)
+ \mathrm{div}(n\nabla S \otimes \nabla S)
+ n\nabla V
= \upsilon^2 \partial_t (n \partial_t S \nabla S),\\[2mm]
-\Delta V = n-b(x) .
\end{cases}
\end{equation}

\medskip
\noindent\textbf{Formal non-relativistic limit.}
When relativistic effects are negligible, corresponding to $\upsilon \to 0$, the relativistic quantum hydrodynamic system \eqref{J1} formally reduces to the quantum Euler--Poisson system
\begin{equation}\label{EP-quantum}
~~~~~\begin{cases}
\partial_t n + \mathrm{div}(n\nabla S) = 0,\\[2mm]
\partial_t (n\nabla S)
+ \mathrm{div}(n\nabla S \otimes \nabla S)
- \dfrac{\varepsilon^2}{2} n \nabla\!\left(\dfrac{\Delta\sqrt{n}}{\sqrt{n}}\right)
+ n\nabla V = 0,\\[2mm]
-\Delta V = n-b(x) .
\end{cases}
\end{equation}

\medskip
\noindent\textbf{Formal non-relativistic--semiclassical limit.}
When both relativistic and quantum effects can be neglected, corresponding to $\upsilon \to 0$ and $\varepsilon \to 0$, the relativistic quantum hydrodynamic system reduces, at a purely formal level, to the classical Euler--Poisson system
\begin{equation}\label{EP-classical}
\begin{cases}
\partial_t n + \mathrm{div}(n\nabla S) = 0,\\[2mm]
\partial_t (n\nabla S)
+ \mathrm{div}(n\nabla S \otimes \nabla S)
+ n\nabla V = 0,~~~~~~~~~~~~~~~~~\\[2mm]
-\Delta V = n-b(x) .
\end{cases}
\end{equation}

The rest of this paper is organized as follows.
In Section~\ref{ii}, we establish a hydrodynamic reformulation of the self-consistent
Klein--Gordon equation via the Madelung transformation, which leads to the
relativistic quantum hydrodynamic system.
In Section~\ref{iii}, we first reformulate the resulting system as a coupled hyperbolic–elliptic system. We then establish a linearized system and derive the corresponding a priori estimates. Finally, we prove the local-in-time existence and uniqueness of classical solutions to the full nonlinear system.

\vskip 0.2cm
\section{Formal derivation of the relativistic quantum hydrodynamic system}\label{ii}
We begin by recalling the physical motivation underlying the modulation of the
Klein--Gordon field in the non-relativistic regime. Since the rest-mass energy
$mc^2$ induces rapid temporal oscillations, it is convenient to factor out this
high-frequency component. To this end, the Klein--Gordon field $\varphi(x,t)$ is
written in the modulated form
\[
\varphi(x,t)=\phi(x,t)\, e^{-\frac{i}{\hbar}mc^2 t},
\]
where the complex-valued function $\phi(x,t)$ varies slowly in time.
In the energy representation, one may formally write
\[
\phi(x,t)=u_E(x)\, e^{-\frac{i}{\hbar}E' t},
\qquad
E' = E - mc^2 = \sqrt{m^2c^4 + c^2 p^2} - mc^2.
\]

In the non-relativistic limit $c\to\infty$, the kinetic energy satisfies
\[
E' \sim \frac{p^2}{2m},
\]
and hence $E' \ll mc^2$, corresponding to the regime $|p|/m \ll c$. In this sense,
the modulated field $\phi(x,t)$ captures the leading-order non-relativistic
dynamics of the Klein--Gordon field after removing the rapidly oscillating
rest-mass energy component, and is consistent with the Schr\"odinger description.

We now turn to the formal derivation of a relativistic quantum hydrodynamic
system starting from the self-consistent Klein--Gordon--Poisson equation.
Let $n=n(x,t)>0$ denote the particle density and $S=S(x,t)$ the phase function
associated with the wave function in spacetime $\mathbb{R}^{3+1}$. Substituting
the above modulation into the Klein--Gordon equation, one formally finds that $\phi$
satisfies
\begin{equation}\label{Hhh1}
i\hbar\partial_t\phi+\frac{\hbar^2}{2m}\Delta\phi
-V(x,t)\phi
=\frac{\hbar^2}{2mc^2}\partial_t^2\phi .
\end{equation}
To identify the key dimensionless parameters governing the relativistic and
quantum effects, we introduce the scalings
\[
x=L\hat{x},\qquad t=T\hat{t},
\]
where $L$ and $T$ denote the characteristic length and time scales, respectively.
Let $U=L/T$ be the reference velocity and rescale the potential as
$V=mU^2\hat{V}$. Substituting these scalings into \eqref{Hhh1} and omitting the
carets yields
\[
i\varepsilon\partial_t\phi+\frac{1}{2}\varepsilon^2\Delta\phi
-V(x,t)\phi
=\frac{1}{2}\varepsilon^2\upsilon^2\partial_t^2\phi,
\]
where the dimensionless parameters
\[
\upsilon=\frac{U}{c},\qquad
\varepsilon=\frac{\hbar}{mU^2T}
\]
measure the relativistic and quantum effects, respectively.
The resulting dimensionless Klein--Gordon--Poisson system reads
\begin{equation}\label{H2}
\begin{cases}
i\varepsilon\partial_t\phi+\frac{1}{2}\varepsilon^2\Delta\phi
-V(x,t)\phi
=\frac{1}{2}\varepsilon^2\upsilon^2\partial_t^2\phi,\\[2mm]
-\Delta V=|\phi|^2-b(x).
\end{cases}
\end{equation}
Here and in what follows, we work with the dimensionless system \eqref{H2}
and suppress the hats for notational simplicity. In particular, the background charge density $b(x)$ is already dimensionless. The WKB analysis (Wentzel \cite{WKB1}, Kramers \cite{WKB2}, Brillouin \cite{WKB3})
provides an effective framework for connecting microscopic quantum dynamics with
macroscopic hydrodynamic behavior. Motivated by this approach, we introduce the
polar decomposition of the wave function. Assuming that $|\phi|>0$, we write
\[
\phi=\sqrt{n}\exp\!\left(\frac{iS}{\varepsilon}\right),
\]
where $n=|\phi|^2$ denotes the particle density and $S$ is the associated phase
function. Substituting this ansatz into the self-consistent Klein--Gordon system
\eqref{H2} yields the corresponding hydrodynamic formulation.
In terms of $\phi$, one has
\[
n=\phi\bar{\phi},
\qquad
\nabla S=\frac{i\varepsilon}{2|\phi|^{2}}
\bigl(\phi\nabla\bar{\phi}-\bar{\phi}\nabla\phi\bigr),
\qquad
\partial_t S=\frac{i\varepsilon}{2|\phi|^{2}}
\bigl(\phi\bar{\phi}_t-\bar{\phi}\phi_t\bigr).
\]
Substituting the decomposition
$\phi=\sqrt{n}\exp(iS/\varepsilon)$
into the Klein--Gordon equation \eqref{H2}$_1$ and dividing by
$\exp(iS/\varepsilon)$ yield
\begin{align}\label{H3}
\frac{i\varepsilon}{2}\frac{\partial_t n}{\sqrt{n}}
-\sqrt{n}\,\partial_t S
&=-\frac{\varepsilon^{2}}{2}
\left(
\Delta\sqrt{n}
+\frac{2i}{\varepsilon}\nabla\sqrt{n}\cdot\nabla S
+\frac{i}{\varepsilon}\sqrt{n}\,\Delta S
-\frac{\sqrt{n}}{\varepsilon^{2}}|\nabla S|^{2}
\right) \notag\\
&\quad
+\sqrt{n}V
+\frac{\varepsilon^2\upsilon^2}{2}\partial_{tt}\sqrt{n}
+i\varepsilon\upsilon^2 \partial_t S\partial_t\sqrt{n}
+\frac{i\varepsilon\upsilon^2}{2}\partial_{tt} S\sqrt{n}
-\frac{\upsilon^2}{2}(\partial_t S)^2\sqrt{n}.
\end{align}
Taking the imaginary part of \eqref{H3} gives
\[
\partial_t n
=-2\sqrt{n}\,\nabla\sqrt{n}\cdot\nabla S
-n\Delta S
+\upsilon^2\partial_t(n\partial_t S)
=-\operatorname{div}(n\nabla S)
+\upsilon^2\partial_t(n\partial_t S),
\]
which coincides with the first equation in \eqref{J1}.
Dividing the real part of \eqref{H3} by $\sqrt n$, differentiating the resulting
equation with respect to $x$, and multiplying by $n$, we obtain, by using the
continuity equation in \eqref{J1},
\begin{align}\label{331sys}
\partial_t(n\nabla S)
+\operatorname{div}\!\left(\frac{n\nabla S\otimes n\nabla S}{n}\right)
&-\frac{\varepsilon^{2}}{2}\,
n\nabla\left(\frac{\Delta\sqrt n}{\sqrt n}\right)
+n\nabla V\notag\\
&=\frac{\upsilon^{2}}{2}\partial_t\bigl(n\partial_t S\nabla S\bigr)
-\frac{\varepsilon^{2}\upsilon^{2}}{2}\,
n\nabla\left(\frac{\partial_{tt}\sqrt n}{\sqrt n}\right).
\end{align}
Equation \eqref{331sys} corresponds to the second equation in \eqref{J1}. Together with the Poisson equation for $V$, this yields the relativistic quantum hydrodynamic system \eqref{J1}. Conversely, let $(n,S)$ be a sufficiently smooth solution of \eqref{J1} with
$n>0$, and define $\phi$ by the same polar decomposition,
\[
\phi=\sqrt n\,\exp\!\left(\frac{iS}{\varepsilon}\right).
\]
A direct computation shows that
\begin{align}
\mathrm{i}\varepsilon \partial_t \phi
+ \frac{\varepsilon^2}{2}\Delta \phi
&= \mathrm{e}^{\mathrm{i}S/\varepsilon}
\Bigg[
\frac{\mathrm{i}\varepsilon}{2\sqrt n}
\Big(
\partial_t n+\operatorname{div}(n\nabla S)
\Big) \notag\\
&\qquad\qquad
+ \sqrt n\Big(
-\partial_t S
- \frac12|\nabla S|^2
+ \frac{\varepsilon^2}{2}
\frac{\Delta\sqrt n}{\sqrt n}
\Big)
\Bigg].
\label{eq:alg}
\end{align}

\medskip
\noindent
We now make use of the Klein--Gordon equation satisfied by $\phi$. More precisely, taking the real and imaginary parts of the Klein--Gordon equation under the Madelung transformation and using \eqref{J1}, we obtain
\begin{align}
\mathrm{i}\varepsilon \partial_t \phi
+ \frac{\varepsilon^2}{2}\Delta \phi
= \sqrt n\,\mathrm e^{\mathrm i S/\varepsilon}
\left(V+\frac12\upsilon^2\varepsilon^2\frac{\partial_{tt}\phi}{\phi}\right)= V\phi+\frac12\upsilon^2\varepsilon^2\partial_{tt}\phi.
\end{align}
 Hence, $\phi$ satisfies the self-consistent
Klein--Gordon system \eqref{H2}. This establishes the formal equivalence between
the Klein--Gordon--Poisson system and its relativistic quantum hydrodynamic
formulation.
The quantum correction term may be interpreted either as a quantum self-potential
(Bohm potential) or, equivalently, as a quantum stress contribution. Indeed,
\[
\frac{\varepsilon^2}{2}\,
n\nabla\!\left(\frac{\Delta\sqrt n}{\sqrt n}\right)
=\frac{\varepsilon^2}{4}\,
\operatorname{div}\!\bigl(n(\nabla\otimes\nabla)\log n\bigr),
\]
where
\[
P=\frac{\varepsilon^2}{4}\,n(\nabla\otimes\nabla)\log n
\]
denotes a non-diagonal quantum stress tensor. The relativistic correction term
accounts for the effect of relativistic particle motion and can be written as
\[
\frac{\varepsilon^{2}\upsilon^{2}}{2}\,
n\nabla\!\left(\frac{\partial_{tt}\sqrt n}{\sqrt n}\right)
=\frac{\varepsilon^{2}\upsilon^{2}}{4}\,
\partial_t\!\left(n\nabla(\partial_t\log n)\right).\]

\begin{remark}
The relativistic quantum hydrodynamic system studied in this paper can be viewed,
at a formal level, as a hydrodynamic reformulation of the self-consistent
Klein--Gordon--Poisson system.
This reformulation helps to clarify the connections between the relativistic
quantum hydrodynamic system and the Euler--Poisson equations under various
parameter regimes.
A rigorous justification of the equivalence between the Klein--Gordon system and
its hydrodynamic formulation, as well as the associated singular limits, is left
for future work.
\end{remark}
\begin{remark}
Assume that the initial data for the Klein--Gordon Cauchy problem satisfy
\[
\phi(x,0)=\phi_0\in H^1(\mathbb{R}^3),\qquad
\partial_t\phi(x,0)=\phi_1\in L^2(\mathbb{R}^3).
\]
On the set where $n_0=|\phi_0|^2>0$ (so that the phase is well-defined), the
corresponding initial data for the relativistic quantum hydrodynamic system
\eqref{J1} can be defined formally by
\[n_0=|\phi_0|^2,\qquad
\nabla S_0=\varepsilon\,\frac{1}{|\phi_0|^2}\,
\operatorname{Im}\bigl(\overline{\phi_0}\,\nabla\phi_0\bigr),\]
\[n_1=\partial_t n(x,0)=\overline{\phi_0}\phi_1+\phi_0\overline{\phi_1},\qquad
S_1=\partial_t S(x,0)
=\varepsilon\,\frac{1}{|\phi_0|^2}\,
\operatorname{Im}\bigl(\overline{\phi_0}\,\phi_1\bigr).\]
Here $S_0$ is determined (locally) up to an additive constant. Conversely, given initial data $(n_0,S_0,n_1,S_1)$ with $n_0>0$, the initial data
for the Klein--Gordon system can be recovered by
\[
\phi_0=\sqrt{n_0}\exp\!\left(\frac{iS_0}{\varepsilon}\right),\qquad
\phi_1=\left(\frac{n_1}{2\sqrt{n_0}}+
\frac{i\sqrt{n_0}S_1}{\varepsilon}\right)
\exp\!\left(\frac{iS_0}{\varepsilon}\right).
\]
\end{remark}

\section{Local existence of a classical solution}\label{iii}
In this section, we establish the local-in-time existence and uniqueness of
classical solutions to the Cauchy problem for the relativistic quantum
hydrodynamic system on a finite time interval.
We consider general smooth initial data, under the assumption that the initial
density is a small perturbation of a positive constant state $\bar n$. More precisely, we study the Cauchy problem associated with system \eqref{J1}
subject to the initial conditions
\begin{align}\label{initial}
(n,S)(x,0) = (n_0,S_0),~(n_t,\partial_t S)(x,0) = (n_1,S_1).
\end{align}
Here, $\bar{n} > 0$ denotes a constant background state.
Throughout this section, we assume that the background charge density
$b(x)$ appearing in the Poisson equation is given by $b(x)=\bar{n}$, and that
the density satisfies the far-field decay condition
\begin{align}\label{far}
\bigl|\sqrt{n(x)}-\sqrt{\bar{n}}\bigr|
\le \frac{C}{|x|^{\alpha}},
\qquad \alpha>3,
\quad \text{as } |x|\to\infty .
\end{align}
The smallness assumption on the density perturbation is introduced to guarantee the strict positivity of $n$, which is essential for the well-posedness of the system due to the presence of the nonlinear quantum
term $\frac{\varepsilon^2}{2}\, n \nabla\!\left(\frac{\Delta \sqrt{n}}{\sqrt{n}}\right)$
and the relativistic correction term $\frac{\varepsilon^{2}\upsilon^{2}}{2}\,
n \nabla\!\left(\frac{\partial_{tt}\sqrt{n}}{\sqrt{n}}\right)$.
Based on a careful analysis of the relativistic quantum hydrodynamic system \eqref{J1}, we are able to establish the following local existence and
uniqueness theorem:
\begin{theorem}[Local existence and uniqueness of classical solutions to the reformulated RQHD system] \label{thm:local}
Assume that there exists $0<\delta\ll \bar n$ such that the initial density satisfies
\[
\bigl|\sqrt{n_0(x)}-\sqrt{\bar n}\bigr|<\sqrt{\delta}
\qquad \text{for all } x\in\mathbb{R}^3,
\]
and the initial data satisfy
\[
(n_0,n_1)\in H^4(\mathbb{R}^3)\times H^3(\mathbb{R}^3),\qquad
(S_0,S_1)\in H^4(\mathbb{R}^3)\times H^3(\mathbb{R}^3).
\]
Then the Cauchy problem \eqref{J1} admits a unique classical solution
$(n,S,V)$ on a short time interval $[0,T^*]$, where $T^*>0$ depends only on
the size of the initial norms. Moreover, the solution satisfies
\[
\sqrt{n(x,t)}\ge \sqrt{\bar n}-\frac{\sqrt{\delta}}{2}
\qquad \text{for all } x\in\mathbb{R}^3,\; 0\le t\le T^*,
\]
and hence
\[
n(x,t)\ge \left(\sqrt{\bar n}-\frac{\sqrt{\delta}}{2}\right)^2>0.
\]
The following regularity properties hold:
$$
\begin{aligned}
\begin{cases}
\sqrt{n}-\sqrt{\bar{n}} \in C\left(\left[0, T^{*}\right]; L^2\left(\mathbb{R}^3\right)\right),\nabla \sqrt{n} \in C\left(\left[0, T^{*}\right]; H^3\left(\mathbb{R}^3\right)\right),\\
\sqrt{n}_t \in C\left(\left[0, T^{*}\right]; H^3\left(\mathbb{R}^3\right)\right),\sqrt{n}_{tt} \in C\left(\left[0, T^{*}\right]; H^2\left(\mathbb{R}^3\right)\right),\\
S \in C\left(\left[0, T^{*}\right]; H^4\left(\mathbb{R}^3\right)\right),\partial_t S \in C\left(\left[0, T^{*}\right]; H^3\left(\mathbb{R}^3\right)\right),\\
\partial_{tt} S \in C\left(\left[0, T^{*}\right]; H^2\left(\mathbb{R}^3\right)\right),\nabla V \in C\left(\left[0, T^{*}\right]; {H}^5\left(\mathbb{R}^3\right)\right).\\
\end{cases}
\end{aligned}
$$
\end{theorem}

\subsection{A hyperbolic--elliptic formulation of the RQHD system}  

To establish the local classical existence for the relativistic quantum hydrodynamics system, we begin by rewriting the system in a more convenient form. For simplicity, we assume that the parameters $\varepsilon=1$ and $\upsilon=1$ are fixed constants; this normalization is adopted only for notational convenience and does not affect the local-in-time well-posedness analysis.

Throughout this section, we work under the condition that the density remains strictly positive, so that $\sqrt{n}$ is well defined. This condition will be justified a posteriori by the lower bound estimate established later.  We introduce the density perturbation variable $\sqrt{n} = \sqrt{\bar n} + \omega$,
so that the relativistic quantum hydrodynamic system \eqref{J1} can be reformulated in terms of the variables $(\omega,S,V)$. To simplify $\eqref{J1}_2$, we first rewrite the divergence term as follows,
\begin{align}
\operatorname{div}(n\nabla S \otimes \nabla S)
&= \frac12\, n \nabla |\nabla S|^2 
   + \operatorname{div}(n\nabla S)\, \nabla S  \\
&= n\,(\nabla\otimes\nabla S)\,\nabla S
   + \operatorname{div}(n\nabla S)\,\nabla S. \notag
\end{align}
Next, we multiply equation $\eqref{J1}_1$ by $\nabla S$ and formally substitute the resulting expression into $\eqref{J1}_2$. This yields the following evolution equation:
\begin{align}
\partial_t S+\frac{1}{2}|\nabla S|^2+V(x,t)-\frac{1}{2}\frac{\Delta \sqrt{n}}{\sqrt{n}}=\frac{1}{2}(\partial_t S)^2-\frac{1}{2}\frac{{\sqrt{n}}_{tt}}{\sqrt{n}}.
\end{align}
Using the density perturbation variable $\omega$, we derive the following coupled hyperbolic--elliptic system for the variables $(\omega,S,V)$. In particular, $\omega$ and $V$ satisfy
\begin{align}\label{omega-linear}
\begin{cases}
\partial_{tt} \omega-\Delta \omega + A_{1}(\omega,V)
= A_{2}(\omega,\partial_t S,\nabla S)
+ A_{3}(\partial_t S,\nabla S),\\[2mm]
-\Delta V = A_{4}(\omega)-\bar n ,
\end{cases}
\end{align}
where the nonlinear terms $A_{1}$–$A_{4}$ are defined by
\begin{align*}
&A_{1}(\omega,V)=2V(\omega+\sqrt{\bar n}),~~A_{2}(\omega,\partial_t S,\nabla S)
=\omega\!\left((\partial_t S)^2-2\partial_t S-|\nabla S|^2\right),\\
&A_{3}(\partial_t S,\nabla S)
=\sqrt{\bar n}\!\left((\partial_t S)^2-2\partial_t S-|\nabla S|^2\right),A_{4}(\omega)=(\sqrt{\bar n}+\omega)^2 .
\end{align*}
By rewriting equation $\eqref{J1}_1$, we further obtain the following hyperbolic equation for the phase function $S$:
\begin{align}\label{S-wave}
\partial_{tt} S-\Delta S
= B_{1}(\omega,\partial_t \omega)
+ B_{2}(\omega,\nabla\omega,\nabla S)
+ B_{3}(\omega,\partial_t \omega,\partial_t S),
\end{align}
where
\begin{align*}
&B_{1}(\omega,\partial_t \omega)
=\frac{2}{\omega+\sqrt{\bar n}}\,\partial_t \omega,~~B_{2}(\omega,\nabla\omega,\nabla S)
=\frac{2}{\omega+\sqrt{\bar n}}\,\nabla\omega\cdot\nabla S,\\
&B_{3}(\omega,\partial_t \omega,\partial_t S)
=\frac{2}{\omega+\sqrt{\bar n}}\,\partial_t\omega\,\partial_t S .
\end{align*}
The initial conditions for the above equations are
\begin{align}
&(\omega,\omega_t)(x,0)=(\omega_0,\omega_1)(x)
:=\Bigl(\sqrt{n_0(x)}-\sqrt{\bar n},\,\frac{n_1(x)}{2\sqrt{n_0(x)}}\Bigr),
\quad x\in\mathbb{R}^3,\\
&\omega(x,t)\to 0 \quad \text{as } |x|\to\infty,~~(S,S_t)(x,0)=(S_0,S_1)(x),~~x\in\mathbb{R}^3.
\end{align}

\begin{remark}
The local existence theory for relativistic quantum fluid systems is highly nontrivial, primarily due to the coexistence of hyperbolic and elliptic structures and the strong coupling between nonlinear relativistic and quantum effects. The quantum terms involve higher–order spatial derivatives, while the relativistic nonlinearities contain higher–order time derivatives. Together they give rise to third–order nonlinear differential operators, which necessitate working in a regime of strictly positive densities and sufficiently regular solutions. Moreover, the maximum principle is not available for obtaining a priori bounds on the density, and it is not apparent how to control the higher–order regularity of the density directly from the equations. For these reasons, we restrict our attention to the short–time existence of classical solutions in a neighborhood of strictly positive densities, employing alternative analytical methods.
\end{remark}

\subsection{A priori estimates}In this subsection, we derive a priori energy estimates for the system.
\begin{definition}[Iteration set]\label{inter}
To construct local-in-time solutions, we introduce the following iteration set.
For constants $M>0$ and $T>0$, define
\[ \mathcal{J}_{M,T} := \Big\{ (\omega,S)\in C([0,T];H^4(\mathbb{R}^3)) \times C([0,T];H^4(\mathbb{R}^3)) \,:\, E(t)\le M \ \text{for all } t\in(0,T] \Big\}, \]
where the energy functional $E(t)$ is defined by
\[
E(t)
:= \|\partial_t \omega(t)\|_{H^3(\mathbb{R}^3)}^2
 + \|\omega(t)\|_{H^4(\mathbb{R}^3)}^2
 + \|\partial_t S(t)\|_{H^3(\mathbb{R}^3)}^2
 + \|S(t)\|_{H^4(\mathbb{R}^3)}^2 .
\]
The precise choice of the time $T>0$ will be specified later so that the
iteration mapping is well defined and becomes a contraction on
$\mathcal{J}_{M,T}$.
\end{definition}

Fix $M>0$ and $T>0$. Let $(\hat{\omega},\hat{S})\in \mathcal{J}_{M,T}$ be given.
Throughout the iteration step, we regard $(\hat{\omega},\hat{S})$ as known functions and solve for $(\omega,S,V)$ from the following linearized system, where the coefficients depending on $(\hat{\omega},\hat{S})$ are frozen.
\begin{align}\label{vKG}
\begin{cases}
\partial_{tt} \omega-\Delta \omega +2V(\omega+\sqrt{\bar{n}})
=(\omega+\sqrt{\bar{n}})\left(  {\partial_{t} \hat{S}}^2-2 \partial_{t} \hat{S}-|\nabla \hat{S}|^2\right),\\
\partial_{tt} S-\Delta S
=\dfrac{2}{(\hat{\omega}+\sqrt{\bar{n}})}
\left(\partial_t \hat{\omega}+\nabla\hat{\omega}\cdot\nabla S-\partial_t \hat{\omega}\partial_t S \right), \\
-\Delta V=(\omega+\sqrt{\bar{n}})^2-\bar{n}.
\end{cases}
\end{align}
We will derive a priori energy estimates for solutions to this linearized
system in the following subsections.

\begin{lemma}[$H^1$ energy estimate]\label{lem:energy-estimate}
Let $(\hat{\omega},\hat{S})\in \mathcal{J}_{M,T}$ be fixed.
Suppose that $(\omega,V,S)\in [C^2(\mathbb{R}^3)]^3$ solves the Cauchy problem
\eqref{vKG}. 
Under the positivity assumption imposed in the previous section, there exists a constant $C=C(\delta)>0$ such that
the following energy estimate holds for all $t\in(0,T]$:
\begin{align}
&\|\partial_t \omega(t)\|_{L^2(\mathbb{R}^3)}^2
 + \|\omega(t)\|_{H^1(\mathbb{R}^3)}^2
 + \|\partial_t S(t)\|_{L^2(\mathbb{R}^3)}^2\notag\\
 &~~~~~~~~~~~~~~~~~~~~~~~+ \|S(t)\|_{H^1(\mathbb{R}^3)}^2
 + \|\nabla V(t)\|_{L^2(\mathbb{R}^3)}^2
\le
E(0)\exp\!\left(\int_0^t C\,\hat{E}(s)\, ds\right),
\end{align}
where
\begin{equation}
\hat{E}(t)
:=
\|\partial_t \hat{\omega}(t)\|_{H^3(\mathbb{R}^3)}^2
+
\|\hat{\omega}(t)\|_{H^4(\mathbb{R}^3)}^2
+
\|\partial_{t} \hat{S}(t)\|_{H^3(\mathbb{R}^3)}^2
+
\|\hat{S}(t)\|_{H^4(\mathbb{R}^3)}^2 .
\end{equation}
\end{lemma}

\begin{proof}
Before proceeding, we recall that the pair $(\hat{\omega},\hat{S})$ belongs to the iteration set $\mathcal{J}_{M,T}$ defined in Definition~\ref{inter}.
We define the $H^1$--level energy functional
\begin{equation}\label{def:E1}
\mathcal E_1(t)
:=
\frac12\Big(
\|\partial_t \omega\|^2_{L^2(\mathbb{R}^3)}
+\|\omega\|^2_{H^1(\mathbb{R}^3)}
+\|\nabla V\|^2_{L^2(\mathbb{R}^3)}
+\|\partial_t S\|^2_{L^2(\mathbb{R}^3)}
+\|S\|^2_{H^1(\mathbb{R}^3)}
\Big).
\end{equation}
For the $\eqref{vKG}_1$ equation, multiplying by $\partial_t \omega$ and using $\eqref{vKG}_3$, then integrating over $\mathbb{R}^3$ yield
\begin{align}\label{oemga1}
\int_{\mathbb{R}^3} \partial_{tt} \omega\partial_t \omega-\Delta \omega \partial_t \omega-V\Delta V_t  
 \mathrm{d} x=\int_{\mathbb{R}^3}\partial_t \omega(\omega+\sqrt{\bar{n}})\left(  {\partial_{t} \hat{S}}^2-2 \partial_{t} \hat{S}-|\nabla \hat{S}|^2\right)\mathrm{d} x.
\end{align}Here we have used the time derivative of the Poisson equation to replace $\Delta V_t$ by terms involving $\omega_t$, so that the elliptic variable $V$ can be absorbed into the hyperbolic energy structure.
Similarly, testing the $\eqref{vKG}_2$ with $\partial_t S$ gives
\begin{align}\label{s1}
~~~~\int_{\mathbb{R}^3} \partial_{tt} S\partial_t S-\Delta S \partial_t S
 \mathrm{d} x=\int_{\mathbb{R}^3}\frac{2}{(\hat{\omega}+\sqrt{\bar{n}})  }\left(\partial_t S\hat{\partial_t \omega}+\partial_t S\nabla\hat{\omega}\cdot\nabla S-\partial_t S\partial_t \hat{\omega}\partial_t S \right)\mathrm{d} x,
\end{align}
Combining \eqref{oemga1} and \eqref{s1}, we derive the following differential energy inequality,
\begin{align}\label{E1}
\frac{1}{2}\frac{d}{dt}\mathcal E_1(t)\leq& C\left(\int_{\mathbb{R}^3}|\partial_t \omega\omega \partial_{t} \hat{S}^2|\mathrm{d}x+\int_{\mathbb{R}^3}|2\partial_t \omega\omega \partial_{t} \hat{S}|\mathrm{d}x+\int_{\mathbb{R}^3}|2\partial_t \omega\omega |\nabla \hat{S}|^2|\mathrm{d}x\right)\notag\\
&+C\sqrt{\bar{n}}\left(\int_{\mathbb{R}^3}|\partial_t \omega \partial_{t} \hat{S}^2|\mathrm{d}x+\int_{\mathbb{R}^3}|2\partial_t \omega \partial_{t} \hat{S}|\mathrm{d}x+\int_{\mathbb{R}^3}|2\partial_t \omega|\nabla \hat{S}|^2|\mathrm{d}x\right)\notag~~~~~~\\
&+C\left(\int_{\mathbb{R}^3}\frac{2}{\hat{\omega}+\sqrt{\bar{n}}}|\partial_t S\hat{\partial_t \omega}|+\frac{2}{\hat{\omega}+\sqrt{\bar{n}}}|\partial_t S\nabla\hat{\omega}\cdot\nabla S|+\frac{2}{\hat{\omega}+\sqrt{\bar{n}}}|\partial_t S\partial_t \hat{\omega} \partial_t S|\right).
\end{align}
By the chain rule in $L^2(\mathbb R^3)$, we obtain
 to $\partial_t S$ and $\partial_t \omega$, we further obtain the estimates
\begin{align}\label{E1}
\frac{d}{dt}\|\omega\|_{L^2\left(\mathbb{R}^3\right)} \leq C\left\|\partial_t \omega\right\|_{L^2\left(\mathbb{R}^3\right)}~and
~\frac{d}{dt}\|S\|_{L^2\left(\mathbb{R}^3\right)} \leq C\left\|\partial_t S\right\|_{L^2\left(\mathbb{R}^3\right)}.
\end{align}
By combining estimates \eqref{E1}, we arrive at
\begin{align}
\frac{1}{2}\frac{d}{dt}\mathcal E_1(t)\leq&C(\|\partial_t \omega\|_{L^2(\mathbb{R}^3)}+\|\partial_t S\|_{L^2(\mathbb{R}^3)})+C\bigg(\int_{\mathbb{R}^3}|\partial_t \omega\omega \partial_{t} \hat{S}^2|\mathrm{d}x+\int_{\mathbb{R}^3}|2\partial_t \omega\omega \partial_{t} \hat{S}|\mathrm{d}x\notag\\
&+\int_{\mathbb{R}^3}|2\partial_t \omega\omega |\nabla \hat{S}|^2|\mathrm{d}x\bigg)+C\sqrt{\bar{n}}\bigg(\int_{\mathbb{R}^3}|\partial_t \omega \partial_{t} \hat{S}^2|\mathrm{d}x+\int_{\mathbb{R}^3}|2\partial_t \omega \partial_{t} \hat{S}|\mathrm{d}x\notag\\
&+\int_{\mathbb{R}^3}|2\partial_t \omega|\nabla \hat{S}|^2|\mathrm{d}x\bigg)
+C\bigg(\int_{\mathbb{R}^3}\frac{2}{\hat{\omega}+\sqrt{\bar{n}}}|\partial_t S\hat{\partial_t \omega}|+\frac{2}{\hat{\omega}+\sqrt{\bar{n}}}|\partial_t S\nabla\hat{\omega}\cdot\nabla S|\notag\\
&+\frac{2}{\hat{\omega}
+\sqrt{\bar{n}}}|\partial_t S\partial_t \hat{\omega} \partial_t S|\bigg)\leq C\sum_{1}^{9} I_i.
\end{align}
By Hölder's inequality together with the Sobolev embedding
\(H^2(\mathbb{R}^3)\hookrightarrow L^\infty(\mathbb{R}^3)\) and
\(H^3(\mathbb{R}^3)\hookrightarrow W^{1,\infty}(\mathbb{R}^3)\),
we estimate the interaction terms \(I_1,\dots,I_9\) as follows.
\begin{align*}
&I_1\leq\int_{\mathbb{R}^3}|\partial_t \omega\omega \partial_{t} \hat{S}^2|\,\mathrm{d}x
\leq C\|\partial_{t} \hat{S}\|^2_{L^{\infty}(\mathbb{R}^3)}
\|\omega\|_{L^2(\mathbb{R}^3)}\|\partial_t \omega\|_{L^2(\mathbb{R}^3)},\\
&I_2\leq\int_{\mathbb{R}^3}|\partial_t \omega\omega \partial_{t} \hat{S}|\,\mathrm{d}x
\leq C\|\partial_{t} \hat{S}\|_{H^2(\mathbb{R}^3)}
\|\omega\|_{L^2(\mathbb{R}^3)}\|\partial_t \omega\|_{L^2(\mathbb{R}^3)},\\
&I_3\leq\int_{\mathbb{R}^3}|\partial_t \omega\omega|\nabla \hat{S}|^2|\,\mathrm{d}x
\leq C\|\hat{S}\|^2_{H^3(\mathbb{R}^3)}
\|\omega\|_{L^2(\mathbb{R}^3)}\|\partial_t \omega\|_{L^2(\mathbb{R}^3)},\\
&I_4\leq\int_{\mathbb{R}^3}|\partial_t \omega\partial_{t} \hat{S}^2|\,\mathrm{d}x
\leq C\|\partial_{t} \hat{S}\|_{H^2(\mathbb{R}^3)}
\|\partial_{t} \hat{S}\|_{L^2(\mathbb{R}^3)}\|\partial_t \omega\|_{L^2(\mathbb{R}^3)},\\
&I_5\leq\int_{\mathbb{R}^3}|\partial_t \omega\partial_{t} \hat{S}|\,\mathrm{d}x
\leq C\|\partial_{t} \hat{S}\|_{L^2(\mathbb{R}^3)}\|\partial_t \omega\|_{L^2(\mathbb{R}^3)},\\
&I_6\leq\int_{\mathbb{R}^3}|\partial_t \omega|\nabla\hat{S}|^2|\,\mathrm{d}x
\leq C\|\hat{S}\|_{H^3(\mathbb{R}^3)}
\|\hat{S}\|_{H^2(\mathbb{R}^3)}\|\partial_t \omega\|_{L^2(\mathbb{R}^3)},\\
&I_7\leq\int_{\mathbb{R}^3}\frac{2}{\hat{\omega}+\sqrt{\bar{n}}}|\partial_t S\partial_t \hat{\omega}|\,\mathrm{d}x
\leq Ca_0\|\partial_t S\|_{L^2(\mathbb{R}^3)}\|\partial_t \hat{\omega}\|_{L^2(\mathbb{R}^3)},\\
&I_8\leq\int_{\mathbb{R}^3}
\frac{2}{\hat{\omega}+\sqrt{\bar{n}}}|\partial_t S\nabla\hat{\omega}\cdot\nabla S|\,\mathrm{d}x
\leq Ca_0\|\hat{\omega}\|_{H^3(\mathbb{R}^3)}
\|\partial_t S\|_{L^2(\mathbb{R}^3)}\|\nabla S\|_{L^2(\mathbb{R}^3)},\\
&I_9\leq\int_{\mathbb{R}^3}
\frac{2}{\hat{\omega}+\sqrt{\bar{n}}}|\partial_t S\partial_t \hat{\omega} \partial_t S|\,\mathrm{d}x
\leq Ca_0\|\partial_t \hat{\omega}\|_{H^2(\mathbb{R}^3)}\|\partial_t S\|^2_{L^2(\mathbb{R}^3)},
\end{align*}
where $a_0:=\|(\hat\omega+\sqrt{\bar n})^{-1}\|_{L^\infty}(\mathbb{R}^3)$,
which is finite due to the positivity assumption. 
By collecting the estimates of $I_1$--$I_9$, we obtain 
\begin{align}\label{eq:E1-Young}
\frac{d}{dt}\mathcal E_1(t)\leq \sum_{i=1}^9 I_i
\le
C\,\hat E(t)\,\mathcal E_1(t).
\end{align}
Applying Grönwall's inequality, we obtain the following energy estimate
\begin{equation}\label{eq:E1-Gronwall}
\mathcal{E}_1(t)
\leq \mathcal{E}_1(0) e^{\int_0^t C \hat{E}(s)\, ds}.
\end{equation}

\end{proof}
We now derive second--order energy estimates for the system.
Although the local existence argument is formulated within an $H^1$--based iteration framework,
the coupling term $V\omega$ gives rise to additional structural features that naturally
emerge at higher--order energy levels.

\begin{lemma}[$H^2$ energy estimate]\label{lem:H2-energy-estimate}
Let $(\hat{\omega},\hat{S})\in \mathcal{J}_{M,T}$ be fixed.
Assume that $(\omega,V,S)$ is a classical solution to the Cauchy problem
\eqref{vKG} on $(0,T]$, i.e.
\[
(\omega,V,S)\in [C^2([0,T]\times\mathbb{R}^3)]^3.
\]
Define the $H^2$--level energy by
\[
\frac12\Big(
\|\nabla\partial_t \omega(t)\|_{L^2(\mathbb{R}^3)}^2
+\|\omega(t)\|_{H^2(\mathbb{R}^3)}^2
+\|\nabla \partial_t S(t)\|_{L^2(\mathbb{R}^3)}^2
+\|S(t)\|_{H^2(\mathbb{R}^3)}^2+\|\nabla V\|^2_{H^1(\mathbb{R}^3)}
\Big).
\]
Then there exists a constant $C=C(\delta)>0$ such that, for all $t\in[0,T]$,
\begin{align}\label{eq:H2-energy}
\frac12\Big(
\|\nabla\partial_t \omega(t)\|_{L^2(\mathbb{R}^3)}^2
+\|\omega(t)\|_{H^2(\mathbb{R}^3)}^2
&+\|\nabla \partial_t S(t)\|_{L^2(\mathbb{R}^3)}^2
+\|S(t)\|_{H^2(\mathbb{R}^3)}^2+\|\nabla V\|^2_{H^1(\mathbb{R}^3)}
\Big)
\notag\\&\leq\mathcal{E}_2(t) \leq \left(\left(\mathcal{E}_2(0)^{-1 / 2}+1\right) e^{-C \hat{E} t / 2}-1\right)^{-2}\end{align}
where
\[
\hat{E}(t)
:=
\|\partial_t \hat{\omega}(t)\|_{H^3(\mathbb{R}^3)}^2
+
\|\hat{\omega}(t)\|_{H^4(\mathbb{R}^3)}^2
+
\|\partial_{t} \hat{S}(t)\|_{H^3(\mathbb{R}^3)}^2
+
\|\hat{S}(t)\|_{H^4(\mathbb{R}^3)}^2 .
\]
\end{lemma}

\begin{proof}
Taking the gradient of the $\omega$–equation and testing it against
$\nabla\partial_t \omega$ in $L^2(\mathbb{R}^3)$ yields
\begin{align}
\int_{\mathbb{R}^3} \nabla\partial_{tt} \omega\cdot\nabla\partial_t \omega-\Delta\nabla \omega\cdot \nabla\partial_t \omega-&\nabla(V \omega)\cdot\nabla\partial_t \omega
 \mathrm{d} x\notag\\
&=\int_{\mathbb{R}^3}\nabla\left(\partial_t \omega(\omega+\sqrt{\bar{n}})( {\partial_{t} \hat{S}}^2-2 \partial_{t} \hat{S}-|\nabla \hat{S}|^2)\right)\mathrm{d} x.
\end{align}
We first apply the product rule to obtain
\begin{equation}\label{eq:product-rule}
\nabla (V\omega)=\nabla V\,\omega+V\nabla\omega.
\end{equation}
Consequently,
\begin{equation}\label{eq:expand-Vomega}
\int_{\mathbb{R}^3}\nabla(V\omega)\cdot\nabla\partial_t \omega\,\mathrm{d}x
=
\int_{\mathbb{R}^3}\nabla V\,\omega\cdot\nabla\partial_t \omega\,\mathrm{d}x
+
\int_{\mathbb{R}^3}V\nabla\omega\cdot\nabla\partial_t \omega\,\mathrm{d}x.
\end{equation}
Next, we make use of the Poisson-type equation satisfied by the potential,
\begin{equation}\label{eq:poisson-Vt}
\nabla\Delta V_t
=
-2\nabla\omega\,\partial_t \omega
-2(\omega+\sqrt{\bar n})\nabla\partial_t \omega,
\end{equation}
which allows us to rewrite the first term on the right-hand side of
\eqref{eq:expand-Vomega}.  
By integration by parts, we obtain
\begin{align}\label{eq:key-cancellation}
\int \nabla V \omega \cdot \nabla \partial_t \omega&=-\frac{1}{2} \int \nabla V \cdot \nabla \Delta V_t-\int \nabla V \cdot \nabla \omega \partial_t \omega\notag\\
&=
\int_{\mathbb{R}^3}
\Big(
-2\nabla V\nabla\omega\,\partial_t \omega
-\sqrt{\bar n}\nabla V\nabla\partial_t \omega
\Big)\,\mathrm{d}x. 
\end{align}
Combining \eqref{eq:expand-Vomega}--\eqref{eq:key-cancellation}, we arrive at the
following energy identity:
\begin{equation}\label{eq:energy-identity}
\begin{aligned}
\frac{\mathrm{d}}{\mathrm{d}t}\frac{1}
{2}\int_{\mathbb{R}^3}
\Big(
|\nabla\partial_t \omega|^2
+|\nabla^2\omega|^2
+|\nabla^2 V|^2
\Big)\,\mathrm{d}x
&=
\int_{\mathbb{R}^3}
\nabla\!\left(
\partial_t \omega(\omega+\sqrt{\bar n})
\big(
\hat \partial_t S^{\,2}
-2\hat \partial_t S
-|\nabla\hat S|^2
\big)
\right)\mathrm{d}x
\\
&\quad
-\int_{\mathbb{R}^3}
\Big(
V\nabla\omega\cdot\nabla\partial_t \omega
-2\nabla\omega\,\partial_t \omega
-\sqrt{\bar n}\,\nabla\partial_t \omega
\Big)\,\mathrm{d}x.
\end{aligned}
\end{equation}
The structure of \eqref{eq:poisson-Vt} yields a favorable cancellation of the
terms involving $\nabla\partial_t \omega$, which plays a crucial role in closing the
energy estimate.
Similarly, testing the $\eqref{vKG}_2$ with $\partial_t S$ gives
\begin{align}\label{sH2}
\frac{\mathrm{d}}{\mathrm{d}t}\frac{1}{2} \int_{\mathbb{R}^3} |\nabla \partial_t S|^2+ |\nabla^2 S|^2
 \mathrm{d} x=\int_{\mathbb{R}^3}\nabla\omega_t\cdot\nabla\left(\frac{2}{(\hat{\omega}+\sqrt{\bar{n}})  }\left(\partial_t S\hat{\partial_t \omega}+\partial_t S\nabla\hat{\omega}\cdot\nabla S-\partial_t S\partial_t \hat{\omega}\partial_t S \right)\right)\mathrm{d} x.
\end{align}
Adding \eqref{eq:energy-identity} and \eqref{sH2}, we define the $H^2$--level
energy functional
\[
\mathcal{E}_2(t)
=
\frac12\int_{\mathbb{R}^3}
\Big(
|\nabla\partial_t \omega|^2
+|\nabla^2\omega|^2
+|\nabla^2 V|^2
+|\nabla \partial_t S|^2
+|\nabla S|^2
\Big)\,\mathrm{d}x.
\]
Then $\mathcal{E}_2(t)$ satisfies
\begin{align}\label{eq:E2}
\frac{d}{dt}\mathcal{E}_2(t)
&=
\int_{\mathbb{R}^3}
\nabla\!\Big(
\partial_t \omega(\omega+\sqrt{\bar n})
(\hat \partial_t S^{\,2}-2\hat \partial_t S-|\nabla\hat S|^2)
\Big)\,\mathrm{d} x
\notag\\
&\quad
+\int_{\mathbb{R}^3}
\nabla\!\Big(
\frac{2}{\hat{\omega}+\sqrt{\bar n}}
\big(
\partial_t S\partial_t \hat{\omega}
+\partial_t S\nabla\hat{\omega}\cdot\nabla S
-\partial_t S^2\partial_t \hat{\omega}
\big)
\Big)\,\mathrm{d} x
\notag\\\
&\quad
-\int_{\mathbb{R}^3}
\Big(
V\nabla\omega\cdot\nabla\partial_t \omega
-2\nabla\omega\,\partial_t \omega
-\sqrt{\bar n}\,\nabla\partial_t \omega
\Big)\,\mathrm{d} x.
\end{align}
By a straightforward computation
\begin{align}
    \frac{d}{dt}\mathcal{E}_2(t)\leq &\int_{R^3}| \nabla(\partial_t \omega(\omega+\sqrt{\bar n})
\hat \partial_t S^{\,2})| \mathrm{d} x+\int_{R^3}|\nabla(\partial_t \omega(\omega+\sqrt{\bar n})2\hat \partial_t S)| \mathrm{d} x
\notag\\&~~~~~~+\int_{R^3}|\nabla(\partial_t \omega(\omega+\sqrt{\bar n})|\nabla\hat S|^2)| \mathrm{d} x
+\int_{R^3}|\nabla
(\frac{2}{\hat{\omega}+\sqrt{\bar n}}
\partial_t S\partial_t \hat{\omega})|dx\notag\\
&~~~~~~+ \int_{R^3}|\nabla(
\frac{2}{\hat{\omega}+\sqrt{\bar n}}\partial_t S\nabla\hat{\omega}\cdot\nabla S)| \mathrm{d} x
+ \int_{R^3}|\nabla
(\frac{2}{\hat{\omega}+\sqrt{\bar n}}\partial_t S^2\partial_t \hat{\omega})|dx\notag\notag\\
&~~~~~~+\int_{\mathbb{R}^3}
|V\nabla\omega\cdot\nabla\partial_t \omega|
+|2\nabla V\nabla\omega\,\partial_t \omega|
+|\nabla V \sqrt{\bar n}\,\nabla\partial_t \omega |\mathrm{d} x\leq \sum_1^9 J_i.
\end{align}
We first note that the $H^2$ estimates for $J_1-J_6$ can be derived in
the same manner as in the $H^1$ case by applying one additional spatial derivative and using standard Sobolev and Moser--type estimates. Since no new structural difficulty arises, we omit the details. We therefore focus on the terms involving the potential $V$, which constitute the only genuinely new contributions at the $H^2$ level.
Using H\"older's inequality and the Gagliardo--Nirenberg inequality, we estimate $J_7$--$J_9$ as follows,
\begin{align*}
J_7
&\le \int_{\mathbb{R}^3} |V|\,|\nabla\omega|\,|\nabla\partial_t \omega| \,\mathrm{d}x\le \|V\|_{L^6(\mathbb{R}^3)} 
      \|\nabla\omega\|_{L^3(\mathbb{R}^3)} 
      \|\nabla\partial_t \omega\|_{L^2(\mathbb{R}^3)}\\& \le C\,\|V\|_{L^6(\mathbb{R}^3)}
      \|\nabla\omega\|_{L^2(\mathbb{R}^3)}^{\frac12}
      \|\nabla^2\omega\|_{L^2(\mathbb{R}^3)}^{\frac12}
      \|\nabla\partial_t \omega\|_{L^2(\mathbb{R}^3)},\\
J_8&\le 2 \int_{\mathbb{R}^3} |\nabla V|\,|\nabla\omega|\,|\nabla\partial_t \omega| \,\mathrm{d}x \le 2\,\|\nabla V\|_{L^6(\mathbb{R}^3)} 
        \|\nabla\omega\|_{L^3(\mathbb{R}^3)} 
        \|\nabla\partial_t \omega\|_{L^2(\mathbb{R}^3)}\\
        &\le C\,\|\nabla V\|_{L^6(\mathbb{R}^3)}
      \|\nabla\omega\|_{L^2(\mathbb{R}^3)}^{\frac12}
      \|\nabla^2\omega\|_{L^2(\mathbb{R}^3)}^{\frac12}
      \|\nabla\partial_t \omega\|_{L^2(\mathbb{R}^3)},\\
J_9&\le \sqrt{\bar n} \int_{\mathbb{R}^3} |\nabla V|\,|\nabla\partial_t \omega| \,\mathrm{d}x \le \sqrt{\bar n}\,
      \|\nabla V\|_{L^2(\mathbb{R}^3)} 
      \|\nabla\partial_t \omega\|_{L^2(\mathbb{R}^3)}.
\end{align*}
As a consequence, there exists a positive constant
$C(\delta)$ such that
\begin{align} \frac d {dt} \mathcal{E}_2(t)\leq C\hat{E}\mathcal{E}_2(t)+C\hat{E}\mathcal{E}_2^{\frac 32}(t). \end{align}
Applying a standard comparison argument to the above nonlinear differential inequality yields the explicit bound \eqref{eq:H2-energy}
\end{proof}
The higher--regularity estimates at the $H^3$ and $H^4$ levels can be derived
in the same manner as the $H^1$ and $H^2$ cases by applying higher--order
spatial derivatives and using standard Sobolev embeddings together with
Moser--type product and commutator estimates.
For brevity, we only state the final $H^4$ a priori bound.
\begin{proposition}[A priori $H^4$ estimate]\label{prop:main-H4}
Let $(\hat{\omega},\hat{S})\in \mathcal{J}_{M,T}$ be fixed.
Assume that $(\omega,V,S)$ is a classical solution to the Cauchy problem
\eqref{vKG} on $(0,T]$, i.e.
\[
(\omega,V,S)\in [C^2([0,T]\times\mathbb{R}^3)]^3.
\]
Define the $H^4$--level energy by
\begin{equation}\label{def:EH4}
\mathcal{E}_4(t):=
\frac12\Big(
\|\partial_t \omega(t)\|_{H^3(\mathbb{R}^3)}^2
+\|\omega(t)\|_{H^4(\mathbb{R}^3)}^2
+\|\nabla V(t)\|_{H^3(\mathbb{R}^3)}^2
+\|\partial_t S(t)\|_{H^3(\mathbb{R}^3)}^2
+\|S(t)\|_{H^4(\mathbb{R}^3)}^2
\Big).
\end{equation}
Then there exists a constant $C=C(\delta)>0$ such that, for all $t\in[0,T]$,
\begin{equation}\label{eq:EH4-main}
\mathcal{E}_4(t) \leq \left(\left(\mathcal{E}_4(0)^{-1 / 2}+1\right) e^{-C \hat{E} t / 2}-1\right)^{-2}.
\end{equation}
\end{proposition}

\subsection{Existence of Solutions for the Linear system}
In this section, we establish the existence theory for the linearized system. We begin by analyzing the hyperbolic-elliptic system consisting of the wave equation for $\omega$ and the Poisson equation for $V$:
\begin{align}\label{vlinear}
\begin{cases}
\partial_{tt} \omega-\Delta \omega +2V(\omega+\sqrt{\bar{n}})=(\omega+\sqrt{\bar{n}})\left(  {\partial_{t} \hat{S}}^2-2 \partial_{t} \hat{S}-|\nabla \hat{S}|^2\right),\\
-\Delta V={(\omega+\sqrt{\bar{n}})}^{2}- \bar n,\\
\omega(x,0)=\omega_0,~\partial_t \omega(x,0)=\omega_1, 
\omega(x,t)\rightarrow0,|x|\rightarrow \infty.
\end{cases}
\end{align}
Before proceeding with the proof, we first recall some useful estimates for the Newtonian potential.
\begin{lemma}[Estimates for Newtonian convolution terms]\label{lem:newton-estimates}
Let $K(x)=\frac{1}{4\pi|x|}$ be the Newtonian kernel in $\mathbb R^3$.
Then the following estimates hold.

\begin{enumerate}[(i)]
\item
If $\omega(\cdot,t)\in L^{12/5}(\mathbb R^3)$, then
\[
K*\omega^2 \in L^6(\mathbb R^3),
\qquad
\|K*\omega^2\|_{L^6(\mathbb R^3)}
\le
C\,\|\omega(\cdot,t)\|_{L^{12/5}(\mathbb R^3)}^{2}.
\]


\item
If $\omega(\cdot,t)\in L^2(\mathbb R^3)\cap L^\infty(\mathbb R^3)$, then
\[
K*\omega^2 \in L^\infty(\mathbb R^3),
\qquad
\|K*\omega^2\|_{L^\infty(\mathbb R^3)}
\le C\left(\|\omega(\cdot,t)\|_{L^\infty(\mathbb R^3)}^{2}
+\|\omega(\cdot,t)\|_{L^2(\mathbb R^3)}^{2}\right).
\]
\end{enumerate}
\end{lemma}

\begin{proof}
\noindent
(i)
Assume $\omega(\cdot,t)\in L^{12/5}(\mathbb R^3)$.
Then $\omega^2(\cdot,t)\in L^{6/5}(\mathbb R^3)$ and
\[
\|\omega^2(\cdot,t)\|_{L^{6/5}}
=
\|\omega(\cdot,t)\|_{L^{12/5}}^{2}.
\]
By the Hardy--Littlewood--Sobolev inequality in dimension $3$
(applied to the kernel $|x|^{-1}$), choosing $p=6/5$ and $q=6$ such that
\[
\frac{1}{q}=\frac{1}{p}-\frac{2}{3},
\]
we obtain
\[
\|K*\omega^2(\cdot,t)\|_{L^6(\mathbb R^3)}
\le
C\,\|\omega^2(\cdot,t)\|_{L^{6/5}(\mathbb R^3)}
\le
C\,\|\omega(\cdot,t)\|_{L^{12/5}(\mathbb R^3)}^{2}.
\]

\medskip
\noindent
(ii)
We split the convolution into near-field and far-field parts:
\[
(K*\omega^2)(x)
=\frac{1}{4\pi}\int_{|x-y|<1}\frac{\omega^2(y)}{|x-y|}\,dy
+\frac{1}{4\pi}\int_{|x-y|\ge1}\frac{\omega^2(y)}{|x-y|}\,dy
=:I_1+I_2.
\]
For the near-field term, using $\omega^2(y)\le\|\omega\|_{L^\infty}^2$,
\[
|I_1|
\le
\frac{1}{4\pi}\|\omega\|_{L^\infty(\mathbb{R}^3)}^2
\int_{|z|<1}\frac{1}{|z|}\,dz
=
\frac12\,\|\omega\|_{L^\infty(\mathbb{R}^3)}^2.
\]
For the far-field term, since $|x-y|\ge1$,
\[
|I_2|
\le
\frac{1}{4\pi}\int_{\mathbb R^3}\omega^2(y)\,dy
=
\frac{1}{4\pi}\|\omega\|_{L^2(\mathbb R^3)}^{2}.
\]
This completes the proof.
\end{proof}
This reduction follows the strategy introduced by Guo \cite{Y. Guo},
where the Poisson equation is incorporated at the operator level in the
analysis of the Euler--Poisson system, and the elliptic variable is
recovered only after the existence of the hyperbolic unknowns has been
established.
\begin{lemma}[Existence of classical solutions to the coupled hyperbolic-elliptic system]
\label{lem:existence-vlinear}
Let $\hat S$ satisfy the regularity condition \eqref{inter}.
Assume that the initial data satisfy
\[
\omega_0 \in H^4(\mathbb{R}^3), \qquad \omega_1 \in H^3(\mathbb{R}^3),
\]
and let $\bar n>0$ be given. Then there exists a time $T_0>0$ such that the system \eqref{vlinear} admits a classical solution 
$(\omega,V)$ on $\mathbb{R}^3\times(0,T_0)$, satisfying
\begin{align*}
&\omega \in C([0,T_0];H^4(\mathbb{R}^3))
     \cap C^1([0,T_0];H^3(\mathbb{R}^3)),\\
&V \in C([0,T_0];C^{2}(\mathbb{R}^3)),~\nabla V \in C([0,T_0];H^5(\mathbb{R}^3)),
\end{align*}
and $(\omega,V)$ satisfies the coupled hyperbolic-elliptic system \eqref{vlinear}
on $\mathbb{R}^3\times[0,T_0]$. Using the equation $\eqref{vlinear}_1$ and the above regularity of $\omega$ and $V$, we obtain \begin{align}\partial_{tt} \omega \in C([0,T_0];H^2(\mathbb{R}^3))\hookrightarrow C([0,T_0];C(\mathbb{R}^3)).\end{align}
Thus $(\omega,V)$ is a classical solution of \eqref{vlinear} on 
$\mathbb{R}^3\times[0,T_0]$.
\end{lemma}
\begin{proof}
We define the potential $V(x,t)$ as a nonlocal functional of $\omega(x,t)$ via the Newtonian potential,
\begin{align}\label{eq:V-formula} 
V(x,t) := -\frac{1}{4\pi} \int_{\mathbb{R}^3} 
\frac{\omega^2(y,t) + 2\sqrt{\bar{n}}\, \omega(y,t)}{|x - y|} \, dy,
\end{align}
which satisfies the corresponding Poisson equation,
\[
-\Delta V(x,t) = (\omega(x,t) + \sqrt{\bar{n}})^2 - \bar{n}.
\]
Thus, the Poisson equation is incorporated at the operator level.
Given that $\omega(x,t) \sim |x|^{-\alpha}$ with $\alpha>3$, the source term decays as 
$|x|^{-2\alpha}+|x|^{-\alpha}$ for large $|x|$. 
By classical estimates for Newtonian potentials with sufficiently fast decaying sources,
$V(x,t)$ satisfies
\[
|V(x,t)| \le C |x|^{-1} \quad \text{as } |x| \to \infty.
\]
Therefore,
\[
V(x,t) \to 0 \quad \text{as } |x| \to \infty.
\]
Substituting \eqref{eq:V-formula} into the wave equation, we obtain the equivalent nonlocal semilinear wave equation,
\begin{align}\label{eq:nonlocal-wave-eq}
    \omega_{tt}(x,t) - \Delta \omega(x,t)
    + 2(\omega(x,t)+\sqrt{\bar{n}})
        \left(
            -\frac{1}{4\pi} \int_{\mathbb{R}^3}
\frac{\omega^2(y,t)+2\sqrt{\bar{n}}\,\omega(y,t)}{|x-y|}\,dy
        \right)=(\omega+\sqrt{\bar{n}})\left(  {\partial_{t} \hat{S}}^2
    -2 \partial_{t} \hat{S}-|\nabla \hat{S}|^2\right).
\end{align}
Conversely, if \( \omega \) is a sufficiently regular solution to
\eqref{eq:nonlocal-wave-eq} and \( V \) is defined by \eqref{eq:V-formula},
then the pair \( (\omega,V) \) solves the original coupled system
\eqref{vlinear} in the classical sense.
Moreover, using the Sobolev embedding
\( H^4(\mathbb{R}^3)\hookrightarrow W^{2,\infty}(\mathbb{R}^3) \),
the algebra property of \( H^4(\mathbb{R}^3) \),
and Lemma~\ref{lem:newton-estimates},
the nonlocal mapping
\[
(\omega+\sqrt{\bar{n}}) \;\mapsto\;
(\omega+\sqrt{\bar{n}})\,\Bigl(K*(\omega^2+2\sqrt{\bar{n}}\omega)\Bigr)
\]
is locally Lipschitz from \( H^4(\mathbb{R}^3) \) to \( H^3(\mathbb{R}^3) \).
The remaining terms on the right-hand side of
\eqref{eq:nonlocal-wave-eq} are smooth and at most quadratic in \( \omega \)
and its derivatives. Therefore, \eqref{eq:nonlocal-wave-eq} is a semilinear wave equation
with locally Lipschitz nonlinearity in the energy space
\( H^4 \times H^3 \). Hence, local well-posedness follows from the standard theory for semilinear wave equations with convolution nonlinearities (see, e.g., Bhimani, EJDE 2021\cite{Bhimani1}).
By elliptic regularity for the Poisson equation on \( \mathbb{R}^3 \) and the Newtonian potential representation, the solution \( V = K * (\omega^2+2\sqrt{\bar{n}}\omega) \) with \( (\omega^2+2\sqrt{\bar{n}}\omega) \in H^4(\mathbb{R}^3) \) satisfies \( \nabla V \in H^5(\mathbb{R}^3) \) and \( V \in W^{2,\infty}(\mathbb{R}^3) \).
\end{proof}

\begin{remark}
We emphasize that the $H^5$ regularity of $\nabla V$ is not required for
the energy estimates nor for the local Lipschitz continuity of the
nonlocal term.
At the $H^4$ energy level it suffices to control $\nabla V$ in
$H^3(\mathbb R^3)$, which already implies $V\in W^{1,\infty}(\mathbb R^3)$
by Sobolev embedding.
The higher regularity $\nabla V\in H^5(\mathbb R^3)$ follows a posteriori from standard elliptic estimates applied to
$-\Delta V=\omega^2+2\sqrt{\bar n}\omega$ and is included only to indicate the classical
smoothness of the solution.
\end{remark}

Similarly, we obtain a classical solution to the hyperbolic equation for $S$. The standard local well-posedness theory for semilinear wave
equations applies (see, e.g., \cite{Evans, sogge}), and the system
\eqref{vlinear} admits a unique local classical solution.
\begin{align}\label{SS}
\begin{cases}
 \partial_{tt} S-\Delta S =\frac{2}{(\hat{\omega}+\sqrt{\bar{n}})  }\left(\hat{\partial_t \omega}+\nabla\hat{\omega}\cdot\nabla S-\partial_t \hat{\omega}\partial_t S \right),\\
 S(x,0)=S_0,\partial_t S(x,0)=S_1.
 \end{cases}
\end{align}
\begin{lemma}[Local existence of classical solutions for the $S$ equation]
Let $\hat{\omega}$ be a given function satisfying \eqref{inter}, and let the
initial data satisfy
\[
(S_0,S_1)\in H^4(\mathbb{R}^3)\times H^3(\mathbb{R}^3).
\]
Then there exists a time $T_1\in[0,T_0]$, depending on
$\|S_0\|_{H^4(\mathbb{R}^3)}$, $\|S_1\|_{H^3(\mathbb{R}^3)}$, and $M$, such that
the Cauchy problem \eqref{SS} admits a  classical solution
\[
S\in C([0,T_1];H^4(\mathbb{R}^3))
   \cap C^1([0,T_1];H^3(\mathbb{R}^3))
   \cap C^2([0,T_1];H^2(\mathbb{R}^3)).
\]
\end{lemma}

\begin{remark}
Let $\omega=\sqrt{n}-\sqrt{\bar n}$. Since
$\omega\in C\bigl([0,T_0];H^4(\mathbb{R}^3)\bigr)
\hookrightarrow C\bigl([0,T_0];C^2(\mathbb{R}^3)\bigr)$,
we have $\|\omega(t)-\omega_0\|_{L^\infty}\to 0$ as $t\to 0$.
If the initial density satisfies
$\sqrt{n_0(x)}=\omega_0(x)+\sqrt{\bar n}\ge \sqrt{\bar n}-\sqrt{\delta}>0$,
then for sufficiently small $T_0>0$,
$\sqrt{n}(x,t)=\omega(x,t)+\sqrt{\bar n}
\ge \sqrt{\bar n}-\frac{\sqrt{\delta}}{2}$.Thus the density remains strictly positive on $[0,T_0]$.
\end{remark}

\subsection{Existence and Uniqueness}
In this subsection, we establish the local-in-time existence and uniqueness of classical solutions to the coupled hyperbolic–elliptic system. The proof follows directly from the energy estimates established in Proposition \ref{lem:energy-estimate} and standard arguments in the energy method.
\begin{lemma}\label{lem:uniform-boundedness}
Let \( (\hat{\omega}, \hat{S}) \in \mathcal{J}_{M,T} \) be fixed. Suppose that \( (\omega, V, S) \in [C^2(\mathbb{R}^3)]^3 \) solves the Cauchy problem for \(\eqref{vKG}\). Then, there exists a constant \( T^* > 0 \), depending on \( M \), such that the following uniform bound holds for all \( t \in [0, T^*] \), where \( T^* \) is chosen as:
\[
T^* = \min\left\{\frac{2}{C(\delta)\hat{E}} \ln \left(\frac{1+E_0^{-1 / 2}}{1+\frac{1}{\sqrt{2}} E_0^{-1 / 2}}\right), 1 \right\}, \quad M = 2E(0).
\]
Additionally, the following energy estimate holds for all \( t \in [0, T^*] \):
\[
E(t) < M,
\]
where \( M \) is a constant depending on the initial energy \( E(0) \) and $\delta$.
\end{lemma}
\begin{proposition}
[Local existence]\label{thm:existence}
Let $M>0$ and let $T^*>0$ be given. Then there exists a time $T^*\in(0,T]$ such that the Cauchy problem \eqref{J1}, \eqref{initial} and \eqref{far} admits at least one solution
$(\omega,S,V)$ on $[0,T^*]\times\mathbb R^3$ satisfying
\begin{align}
&\omega,S\in C([0,T^*];H^4(\mathbb R^3)), \partial_t \omega,\partial_t S\in C([0,T^*];H^3(\mathbb R^3)),\\
&\partial_{tt} \omega,\partial_{tt} S\in C([0,T^*];H^2(\mathbb R^3)),\nabla V\in C([0,T^*];H^5(\mathbb{R}^3)).
\end{align}
Moreover, $(\omega,S)$ belongs to the invariant set
$\mathcal J_{T^*,M}$.
\end{proposition}
\begin{proof}
We recall from the previous analysis that the electrostatic potential
$V(x,t)$ can be eliminated through the Poisson equation, which reduces
the original system to the coupled quasilinear wave equations
\eqref{eq:nonlocal-wave-eq} and \eqref{SS}.
Therefore, it suffices to establish the existence of solutions
$(\omega,S)$ to this reduced system; the existence of $V$ then follows
directly from elliptic regularity.
Let
\[
B:=L^2([0,T^*];L^2(\mathbb R^3))\times L^2([0,T^*];L^2(\mathbb R^3)),
\]
and recall the definition of the closed convex set
$\mathcal J_{T^*,M}\subset B$ introduced in proposition~\ref{prop:main-H4}. The space $B$ is chosen so that $\mathcal J_{T^*,M}$ is compact in $B$
by the Aubin--Lions lemma, thanks to the uniform $H^4$-bounds.
By construction, $\mathcal J_{T^*,M}$ is convex and bounded in $B$.
For any $(\hat\omega,\hat S)\in\mathcal J_{T^*,M}$, we define the mapping
\[
\mathcal T(\hat\omega,\hat S):=(\omega,S),
\]
where $\omega=\mathcal T_1(\hat S)$ is the solution to the
linearized equation \eqref{eq:nonlocal-wave-eq} with coefficient $\hat S$,
and $S=\mathcal T_2(\hat\omega)$ is the solution to the linearized
equation \eqref{SS} with coefficient $\hat\omega$.
Standard energy estimates imply that $\mathcal T$ maps
$\mathcal J_{T^*,M}$ into itself.
We next show that $\mathcal T$ is continuous from $\mathcal J_{T^*,M}$
into $B$.
Let $\{(\hat\omega_k,\hat S_k)\}_{k\ge1}\subset\mathcal J_{T^*,M}$ be a
sequence converging to $(\hat\omega,\hat S)$ in $B$, and denote
\[
(\omega^k,S^k):=\mathcal T(\hat\omega_k,\hat S_k).
\]
By the uniform $H^4$ a proposition \ref{prop:main-H4} bounds in prop, together
with the Aubin--Lions compactness lemma, there exists a subsequence
(still denoted by $k$) such that, for any $\sigma\in(0,1)$,
\begin{align*}
&S^k \to S ~\text{in } C([0,T^*];H^{4-\sigma}(\mathbb R^3)),
S^k_t \to \partial_t S~\text{in } C([0,T^*];H^{3-\sigma}(\mathbb R^3)),\\
&S^k_{tt} \to \partial_{tt} S~\text{in } C([0,T^*];H^{2-\sigma}(\mathbb R^3)),
\omega^k \to \omega~\text{in } C([0,T^*];H^{4-\sigma}(\mathbb R^3)),\\
&\omega^k_t \to \partial_t \omega~\text{in } C([0,T^*];H^{3-\sigma}(\mathbb R^3)),
\omega^k_{tt}\to \partial_{tt} \omega~\text{in } C([0,T^*];H^{2-\sigma}(\mathbb R^3)).
\end{align*}
Passing to the limit in the linearized equations shows that
$(\omega,S)=\mathcal T(\hat\omega,\hat S)$, which proves the continuity of
$\mathcal T$ on $\mathcal J_{T^*,M}$.
Since $\mathcal J_{T^*,M}$ is convex, closed, and compact in $B$, and
$\mathcal T:\mathcal J_{T^*,M}\to\mathcal J_{T^*,M}$ is continuous, the
Schauder fixed point theorem yields the existence of a fixed point
$(\omega,S)\in\mathcal J_{T^*,M}$ such that
\[
(\omega,S)=\mathcal T(\omega,S).
\]
This fixed point is a solution to the original quasilinear system on
$[0,T^*]$ with the stated regularity. Finally, we recover the electrostatic potential through the Newtonian
potential operator. Define
\[
f(\omega):=(\omega+\sqrt{\bar n})^2-\bar n
=\omega^2+2\sqrt{\bar n}\,\omega.
\]
Under this condition, \(V\) is uniquely given by the Newtonian potential
\[
V(x)=-\frac{1}{4\pi}\int_{\mathbb R^3}\frac{f(\omega)(y)}{|x-y|}\,dy.
\]
By homogeneous elliptic regularity, we obtain
\[
\|\nabla V\|_{H^5(\mathbb R^3)}
\le C\,\|f(\omega)\|_{H^4(\mathbb R^3)},
\]
and consequently
\[
\nabla V\in C((0,T^*];H^5(\mathbb R^3)),
\]
which completes the proof.
\end{proof}
\begin{proposition}[Uniqueness]\label{Uniqueness}
Let $T^*>0$ be the existence time given by Proposition~\ref{thm:existence}.
Assume that $(S_1,\omega_1,V_1)$ and $(S_2,\omega_2,V_2)$
are two classical solutions on $[0,T^*]\times\mathbb R^3$
to the system \eqref{omega-linear} and \eqref{S-wave}
with the same initial data at $t=0$.
Then
\[
(S_1,\omega_1,V_1)\equiv (S_2,\omega_2,V_2)
\quad \text{on } [0,T^*]\times\mathbb R^3 .
\]
\end{proposition}

\begin{proof}
Let $(\omega_i,S_i,V_i)$, $i=1,2$, be two classical solutions to \eqref{omega-linear} and \eqref{S-wave}
on $[0,T^*]\times\mathbb R^3$ with the same initial data. Define the differences
\begin{equation}\label{eq:diff-def}
\eta:=S_1-S_2,\qquad \zeta:=\omega_1-\omega_2,\qquad \Phi:=V_1-V_2 .
\end{equation}
Since the initial data coincide, we have
\begin{equation}\label{eq:E0}
\zeta(0)=\zeta_t(0)=\eta(0)=\eta_t(0)=0,
\end{equation}
and thus the initial energy defined below vanishes. Define the  energy
\begin{equation}\label{eq:energy}
 \tilde{E}(t):=\|\zeta_t(t)\|^2_{L^2(\mathbb{R}^3)}+\|\zeta(t)\|_{H^1(\mathbb{R}^3)}^2
+\|\eta_t(t)\|_{L^2(\mathbb{R}^3)}^2+\|\eta(t)\|_{H^1(\mathbb{R}^3)}^2.
\end{equation}
We test the Poisson equation with $\Phi$. Multiply the $\zeta$-equation by $\zeta_t$ and the $\eta$-equation by $\eta_t$ and use integration by parts to get
\begin{equation}\label{eq:basic-energy}
\frac{d}{dt}\Big(\|\zeta_t(t)\|^2_{L^2(\mathbb{R}^3)}+\|\zeta(t)\|_{H^1(\mathbb{R}^3)}^2
+\|\eta(t)\|_{H^1(\mathbb{R}^3)}^2+\|\eta_t(t)\|_{L^2(\mathbb{R}^3)}^2\Big)
+ \|\nabla\Phi\|_2^2
\le \sum_{k=1}^4 g_k,
\end{equation}
where
\begin{align*}
&g_1=\int_{\mathbb{R}^3}\zeta_t\left(-2V_1(\omega_1+\sqrt{\bar{n}})+2V_2(\omega_2+\sqrt{\bar{n}})\right)\mathrm{d}x,\\
&g_2=\int_{\mathbb{R}^3}\Phi\left((\omega_1+\sqrt{\bar{n}})^2-(\omega_2+\sqrt{\bar{n}})^2\right)\mathrm{d} x,\\
&g_3=\int_{\mathbb{R}^3}\zeta_t\left((\omega_1+\sqrt{\bar{n}})\left(  {\partial_{t} S_1}^2-2 \partial_{t} S_1-|\nabla S_1|^2\right)-(\omega_2+\sqrt{\bar{n}})\left(  {\partial_{t} S_2}^2-2 \partial_{t} S_2-|\nabla S_2|^2\right)\right) \mathrm{d} x,\\
&g_4=\int_{\mathbb{R}^3}\eta_t \bigg(\frac{2}{(\omega_1+\sqrt{\bar{n}})}
\left(\partial_t \omega_1+\nabla\omega_1\cdot\nabla S_1-\partial_t\omega_1\partial_t S_1 \right)-\frac{2}{(\omega_2+\sqrt{\bar{n}})}
\left(\partial_t \omega_2+\nabla\omega_2\cdot\nabla S_2-\partial_t\omega_2\partial_t S_2\right)\bigg)\mathrm{d} x.
\end{align*}
We denote by $C_M$ a constant such that, for $t\in[0,T^*]$,
\[
\|\omega_i(t)\|_{W^{1,\infty}(\mathbb{R}^3)}+\|S_i(t)\|_{W^{1,\infty}(\mathbb{R}^3)}+\|V_i(t)\|_{L^\infty(\mathbb{R}^3)}
+\Big\|\frac1{\omega_i(t)+\sqrt{\bar n}}\Big\|_{L^\infty(\mathbb{R}^3)}\le C_M,\qquad i=1,2.
\]
Such a bound holds for classical solutions on $[0,T^*]$.
Hence, by H\"older and Young inequalities and Sobolev embedding
\begin{align*}
g_1&\leq \int_{\mathbb{R}^3}-2(\omega_1+\sqrt{\bar{n}})\Phi\zeta_t \mathrm{d}x+\int_{\mathbb{R}^3}-2\zeta\zeta_tV_2\mathrm{d}x\notag\\
&\leq C(\|2(\omega_1+\sqrt{\bar{n}})\|_{L^3(\mathbb{R}^3)}\|\Phi\|_{L^6(\mathbb{R}^3)}\|\zeta_t\|_{L^2(\mathbb{R}^3)}+\|V_2\|_{L^{\infty}\mathbb(R^3)}\|2\zeta\|_{L^2\mathbb(R^3)}\|\partial_t\zeta\|_{L^2(\mathbb{R}^3)}\notag\\
&\leq C_M\left(\|\Phi\|_{L^6(\mathbb{R}^3)}\|\partial_t \zeta\|_{L^2(\mathbb{R}^3)}+\|2\zeta\|_{L^2\mathbb(R^3)}\|\partial_t \zeta\|_{L^2(\mathbb{R}^3)}\right),\\
g_2&\leq\int_{\mathbb{R}^3}\Phi\zeta(\omega_1+\omega_2+2\sqrt{\bar{n}})\\
&\leq\|\Phi\|_{L^6(\mathbb{R}^3)}\|\zeta\|_{L^2(\mathbb{R}^3)}\|\omega_1+\omega_2+2\sqrt{\bar{n}}\|_{L^3(\mathbb{R}^3)}\leq C_M\|\Phi\|_{L^6(\mathbb{R}^3)}\|\zeta\|_{L^2(\mathbb{R}^3)},\\
g_3&\leq \int_{\mathbb{R}^3}\partial_t\zeta\zeta({\partial_{t} S_1}^2-2 \partial_{t} S_1-|\nabla S_1|^2)\\
&\qquad+\partial_t\zeta(\omega_2+\sqrt{\bar{n}})\left((\partial_tS_1+\partial_tS_2)\partial_t\eta-2\partial_t\eta-(|\nabla S_1|^2+|\nabla S_2|^2)\nabla \eta\right) \mathrm{d} x\\
&\leq C_M\left(\|\partial_t\zeta\|_{L^2\mathbb(R^3)}\|\zeta\|_{L^2\mathbb(R^3)}+\|\partial_t\zeta\|_{L^2\mathbb(R^3)}\|\partial_t \eta\|_{L^2\mathbb(R^3)}+\|\nabla \eta\|_{L^2\mathbb(R^3)}\|\partial_t \zeta\|_{L^2\mathbb(R^3)}\right),\\
g_4&\leq \int_{\mathbb{R}^3}\partial_t \eta \bigg(-\frac{2\zeta}{(\omega_1+\sqrt{\bar{n}})(\omega_2+\sqrt{\bar{n}})}
\left(\partial_t \omega_1+\nabla\omega_1\cdot\nabla S_1-\partial_t\omega_1\partial_t S_1 \right)\bigg)\\
&\qquad-\partial_t \eta\bigg(\frac{2}{(\omega_2+\sqrt{\bar{n}})}
\bigg(\partial_t \zeta+\nabla S_1\nabla\zeta+\nabla \omega_2\nabla \eta -\partial_tS_1\partial_t\zeta -\partial_t \omega_2 \partial_t \eta \bigg)\bigg)\mathrm{d} x\\
&\leq C_M\bigg(\|\partial_t\eta\|_{L^2(\mathbb{R}^3)}\|\zeta\|_{L^2(\mathbb{R}^3)}+\|\partial_t\eta\|_{L^2(\mathbb{R}^3)}\|\partial_t\zeta\|_{L^2(\mathbb{R}^3)}+\|\partial_t\eta\|_{L^2(\mathbb{R}^3)}\|\nabla \eta\|_{L^2(\mathbb{R}^3)}\\
&\qquad+\|\partial_t\eta\|_{L^2(\mathbb{R}^3)}\|\nabla \eta\|_{L^2(\mathbb{R}^3)}+\|\partial_t\eta\|_{L^2(\mathbb{R}^3)}^2 \bigg).
\end{align*}
Choosing $\varepsilon>0$ sufficiently small and inserting the above estimate, we conclude that
\begin{align}
\frac{d}{dt}\tilde{E}(t)+\|\nabla \Phi\|^2_{L^2(\mathbb{R}^3)}\leq C_M\tilde{E}(t).
\end{align}
Since $\tilde E(0)=0$, Gr\"onwall's inequality yields $\tilde E(t)\equiv 0$ on $[0,T^*]$.
Integrating the differential inequality in time we obtain
\[
\int_0^{T^*}\|\nabla\Phi(t)\|_2^2\,dt=0,
\]
hence $\|\nabla\Phi(t)\|_2=0$ for a.e. $t\in[0,T^*]$.
By the continuity of $t\mapsto \nabla\Phi(t)$ in $L^2(\mathbb R^3)$, it follows that $\nabla\Phi\equiv 0$ on $[0,T^*]\times\mathbb R^3$.
Consequently,
\begin{align}
\zeta\equiv 0,\qquad \eta\equiv 0 \quad \text{on }[0,T^*]\times\mathbb R^3.
\end{align}
By the normalization of the Poisson problem, for each $t\in[0,T^*]$,
$V(t,\cdot)$ is the unique solution of the Poisson equation satisfying
$V(t,x)\to0$ as $|x|\to\infty$. This completes the proof of uniqueness.
\end{proof}
Based on the results established above, we now give the proof of the theorem.
\begin{proof}
The local existence of a classical solution $(\omega,S,V)$ on
$[0,T^*]\times\mathbb R^3$ with the stated regularity follows directly from
Proposition~\ref{thm:existence}.
To prove uniqueness, let
$(\omega_1,S_1,V_1)$ and $(\omega_2,S_2,V_2)$ be two classical solutions
corresponding to the same initial data.
Then Proposition~\ref{Uniqueness} implies that
\[
(\omega_1,S_1,V_1)\equiv(\omega_2,S_2,V_2)
\quad \text{on } [0,T^*]\times\mathbb R^3.
\]
This completes the proof.
\end{proof}
\section*{Acknowledgements}
The research of Ben Duan was supported in part by the National Key R\&D Program
of China (No. 2024YFA1013303) and by the National Natural Science Foundation of
China (Grant Nos. 12271205 and 12171498). The research of Bin Guo was supported
in part by the National Key R\&D Program of China (No. 2024YFA1013301). The
research of Jun Li was supported in part by the National Key R\&D Program of
China (No. 2024YFA1013301). The research of Rongrong Yan was supported in part
by the National Key R\&D Program of China (No. 2024YFA1013303).

\section*{Data Availability}
The authors confirm that the data supporting the findings of this study are available within the article.

\section*{Declarations}
\section*{Conflict of interest}  All authors declare that they have no conflict of interest.

\end{document}